


\magnification=\magstep 1
\parindent = 0 pt
\baselineskip = 16 pt
\parskip = \the\baselineskip

\font\AMSBoldBlackboard = msbm10

\def\RR{{\hbox{\AMSBoldBlackboard R}}}
\def\CC{{\hbox{\AMSBoldBlackboard C}}}

\def\QQ{{\hbox{\AMSBoldBlackboard Q}}}

\def\ZZ{{\hbox{\AMSBoldBlackboard Z}}}

\settabs 12\columns
\rightline{math.NT/9812012}
\vskip 1 true in
{\bf \centerline{SPECTRAL ANALYSIS OF THE LOCAL COMMUTATOR OPERATORS}}
\vskip 0.5 true in
\centerline{Jean-Fran\c{c}ois Burnol}
\par
\centerline{December 1998}
\par
The spectral analysis of the (local) conductor operator $H = \log(|q|) + \log(|p|)$ was shown in a previous paper to be given by the Explicit Formula. I give here the spectral analysis of the commutator operator $K = i\ [\log(|p|),\ \log(|q|)]$ (which shares with $H$ the property of complete dilation invariance). Its spectral function is found to be the derivative of the one of $H$. It is then deduced that $K$ anticommutes with the Inversion and is a bounded operator. Higher commutators are related in the same manner to the higher derivatives of the Gamma function, which thus all acquire an operator theoretic and spectral interpretation.\par
\vfill
{\parskip = 0 pt
62 rue Albert Joly\par
F-78000 Versailles\par
France\par}

jf.burnol@dial.oleane.com

\eject

{\bf
TABLE OF CONTENTS\par
\par
On operators that are completely dilation invariant\par
Explicit spectral analysis over the $p$-adics\par
The Gamma function, its first and second derivatives\par
Spectral analysis in the general case\par
Higher derivatives and higher commutators\par
Looking back and making it all (too) easy\par
References\par}
\par
\vfill \eject

{\bf On operators that are completely dilation invariant}

In [1a,b] it was shown that the local contributions to the Explicit Formula of Analytic Number Theory have an operator theoretic and spectral interpretation in terms of the conductor operator $H = \log(|y|_\nu) + \log(|x|_\nu)$. This operator acts as an hermitian (unbounded, but bounded below) operator on the Hilbert space $L^2(K_\nu, dy)$ on the local field $K_\nu$ with respect to the additive measure $dy$. $K_\nu$ is the completion of an algebraic number field at some place (either finite or archimedean). I have used the physicist's notation so that $y$ and $x$ are Fourier dual variables (with respect to the duality $\lambda(xy) = \exp(\pm2\pi i\ Tr(xy))$, with a $+$ at finite places and a $-$ at archimedean places), and the additive Haar measure $dy$ is always chosen to be self-dual. As in this paper only the local case is considered, I will mostly drop the $\nu$ subscripts.

The key feature allowing the spectral analysis of $H$ is that it is completely dilation invariant. Let $A$ be the operator of multiplication with $\log(|y|)$ and $B$ its Fourier conjugate which acts in the $x$-representation by multiplication with $\log(|x|)$. Let $R(u): \varphi(y) \mapsto {|u|}^{-1/2}\,\varphi(y/u)$. Then $R(u)\,A\,R(u)^{-1} = A - \log(|u|)$, $R(u)\,B\,R(u)^{-1} = B + \log(|u|)$ so that both $H$ and $K = i\,[B,\,A]$ commute with the dilation operators $R(u)$.

This complete dilation invariance of an operator $O$ means that it first needs to be transported to $L^2(K_\nu^\times, d^*u)$ for its analysis (the measure $d^*u$ is chosen to assign a volume of $1$ to the units at a finite place and to be push-forwarded to $dt/t$ under $u \mapsto t = |u|$ at an archimedean place). If $O$ is sufficiently manageable, it should have as generalized eigenvectors the multiplicative unitary characters $\chi(u)$, which from the additive side take the shape $\chi(y)\,|y|^{-1/2}$. 

The spectral analysis of $O$ is then expected to be given as an equation
$$O({\chi(y)}^{-1}|y|^{-s}) = O(\chi, s)\ {\chi(y)}^{-1}|y|^{-s}\leqno ({\cal S}_O)$$
for some functions $O(\chi, s)$ on the critical line $Re(s) = 1/2$. The way we understand these equations is as an identity of distributions on the {\it punctured} $\nu$-adics. We will not address the question of the validity of the identity on the complete line for the commutator operator (it works for the conductor operator).

In [1b] I showed that this program works for the conductor operator, with the result that $$H(\chi, s) = {\partial\over\partial s}\log(\Gamma(\chi,s))$$
where $\Gamma(\chi,s)$ is the Tate-Gel'fand-Graev Gamma function ([2a,b,c]), and the identity makes sense for $s$ in the critical strip but acquires its Hilbert space spectral interpretation only on the critical line. An explicit evaluation then shows that $H$ is bounded below (the exact values will be given below). In this paper I prove that a similar method works for the commutator operator $K = i\,[B, A]$ too, with the result
$$K(\chi, s) = - i {\partial^2\over\partial s^2}\log(\Gamma(\chi,s))$$
and this implies, as will be shown, that $K$ is bounded.

This last result obviously should be related to the theory of Landau, Pollak, and Slepian on band- and time-limited functions [3a,b,c] mentioned in Connes's paper [3d] and in the book of Dym-McKean [3e]. I hope to return to this subject on another occasion.

In the same manner that the conductor operator is related to the Riemann Hypothesis, the commutator operator is related to the statistics of the zeros. I hope to return to these subjects, too.

{\bf Explicit spectral analysis over the $p$-adics}

For concreteness I will first elucidate the spectrum of the commutator operator over $\QQ_p$. As it commutes with averaging over the units one can split the discussion between the ramified and the non-ramified spectrum. I showed in [1a] that it vanishes on the ramified part of $L^2(\QQ_p, dy)$, so that it is enough to look at what happens on the invariant part $L^2_{inv}(\QQ_p, dy)$ spanned by functions depending only on the norm $|y|$.

Let's first evaluate $K(\theta) = i\,(BA - AB)(\theta)$ for $\theta$ the characteristic function of the units. As $A(\theta) = 0$ this is just $-i\,A(H(\theta))$. In [1b] we found
$$\eqalign{
(|y| < 1)\quad\ H(\theta)(y) &= -\log(p) \cr
(|y| = 1)\quad\ H(\theta)(y) &= \;0 \cr
(|y| > 1)\quad\ H(\theta)(y) &= -{\log(p)\over|y|} }$$
so 
$$\eqalign{
(|y| < 1)\quad\ K(\theta)(y) &= i\,\log(p) \log(|y|)\cr
(|y| = 1)\quad\ K(\theta)(y) &= \;0 \cr
(|y| > 1)\quad\ K(\theta)(y) &= i\,\log(p){\log(|y|)\over|y|} }$$

An orthonormal basis of $L^2_{inv}(\QQ_p,dy)$ is given by the functions
$$\eta_j(y) = (1 - 1/p)^{-1/2}\ p^{-j/2}\ \theta(p^j y)$$
By dilation invariance one has $K(\theta_j)(y) = K(\theta)(p^jy)$ for $\theta_j(y) = \theta(p^jy)$, so the matrix elements of $K$ in the orthonormal basis $\{\eta_j, j\in\ZZ\}$ depend only on $j-k$ and are
$$\eqalign{
(j\neq k)\quad\ K_{jk} &=\int_{\QQ_p}\overline{K\eta_j(y)}\eta_k(y)\,dy\ = -i\,(\log(p))^2\,(k-j)\,p^{-|k-j|/2} \cr
(j = k)\quad\ K_{jj} &= 0 }$$
We can now use the isometry from $L^2_{inv}(\QQ_p, dy)$ to $L^2(S^1, {d\theta\over2\pi})$ which sends $\eta_i$ to $e_i(z) : z\mapsto z^i$. Under this isometry one deduces from the above results that $K$ is transported to the multiplication operator with the function $$k(z) = i\log(p)^2\ \sum_{k\neq0}k\,p^{-|k|/2}e_k(z)$$
This is $-i\,\log(p)\,z{\partial\over\partial z}h(z)$ where
$$h(z) = -\log(p)\,\left({z\over\sqrt{p}-z}+{\overline{z}\over\sqrt{p}-\overline{z}}\right)$$
is the spectral function for the conductor operator ([1b]). So
$$k(z) = i\,\log(p)^2\,\left({{z\over\sqrt{p}}\over{(1 - {z\over\sqrt{p}})^2}}-{{\overline{z}\over\sqrt{p}}\over{(1 - {\overline{z}\over\sqrt{p}})^2}}\right)$$
A little calculation then completes the proof of the following statement:

{\bf Theorem:} The ramified spectrum of the commutator operator is reduced to $\{0\}$ and has infinite multiplicity, while the invariant spectrum is continuous with support the closed interval $[-{2\log(p)^2\over\sqrt{p}},+{2\log(p)^2\over\sqrt{p}}]$. The commutator operator is hermitian, bounded, and anticommutes with the Inversion.

{\bf The Gamma function, its first and second derivatives}

The basis of our method remains the fundamental Fourier Transform identity ([2a,b,c])
$${\cal F}(\chi(x)|x|^{s-1}) = \Gamma(\chi, s)\chi^{-1}(y)|y|^{-s}$$
Here $\chi$ is a multiplicative (quasi-)character and $\chi(x)|x|^{s-1}$ and $\chi^{-1}(y)|y|^{-s}$ are analytic continuations from a suitable open strip in the variable $s$. We will always take $\chi$ to be a unitary character, then for $s$ in the critical strip this has the direct meaning that
$$\int \psi(x)\,\chi(x)|x|^{s-1} dx = \Gamma(\chi, s) \int \varphi(y)\,\chi^{-1}(y)|y|^{-s} dy\leqno{(I)}$$
for any Schwartz-Bruhat function $\varphi(y)$ on $K_\nu$ and $\psi(x) = \widetilde{\varphi}(x) = \int_{K_\nu}\varphi(y)\lambda(-xy)\,dy$ its Fourier Transform.

The function $\Gamma(\chi, s)$ of the character $\chi$ is meromorphic with respect to the complex parameter $s$ and of course depends on them only through their combination $\chi(x)|x|^s$. It is analytic and non-vanishing in the critical strip ($\chi$ is unitary), and is periodic with period $2\pi i\over \log(q)$ if $\nu$ is finite, with $q$ the cardinality of the residue field.

In the next chapter, we take the first and then the second derivatives of $(I)$ which give the following identities of distributions on the punctured line (for $0<Re(s)<1$): 
$$H({\chi(y)}^{-1}|y|^{-s}) = H(\chi, s)\ {\chi(y)}^{-1}|y|^{-s}
 = \left({\partial\over\partial s}\log(\Gamma(\chi,s))\right)\ {\chi(y)}^{-1}|y|^{-s}\leqno ({\cal S}_H)$$
$$K({\chi(y)}^{-1}|y|^{-s}) = K(\chi, s)\ {\chi(y)}^{-1}|y|^{-s}
= \left(- i\,{\partial^2\over\partial s^2}\log(\Gamma(\chi,s))\right)\ {\chi(y)}^{-1}|y|^{-s}\leqno ({\cal S}_K)$$
The functions $H(\chi, s)$ and $K(\chi, s)$ have the following properties
$$H(\chi,s) = H(\chi^{-1}, 1-s)$$
$$\overline{H(\chi,s)} = H(\overline{\chi},\overline{s})$$
$$Re(s) = {1\over2}\Rightarrow H(\chi,s) = H(\overline{\chi},\overline{s})\hbox{ and is real}$$
$$K(\chi,s) = - K(\chi^{-1}, 1-s)$$
$$\overline{K(\chi,s)} = - K(\overline{\chi},\overline{s})$$
$$Re(s) = {1\over2}\Rightarrow K(\chi,s) = -K(\overline{\chi},\overline{s})\hbox{ and is real}$$
As the Inversion sends a unitary character to its complex conjugate, this implies:

{\bf Theorem:} The commutator operator anticommutes with the Inversion.

At a finite place the Gamma function of a ramified character was evaluated by Tate [2a]: its logarithmic derivative is a constant so that the next derivative vanishes. In this manner we recover the result of [1a] that $\log(|x|)$ and $\log(|y|)$ commute on the ramified part of the Hilbert space. There only remains the invariant part, which is dual to the circle of non-ramified unitary characters. The spectral functions of $H$ and $K$ are continuous hence bounded. So the commutator operator at a finite place is a bounded hermitian operator, and it remains to check that this works for the real and complex places too.

For the real (resp{.} complex) place there is a canonical way to choose a base point $\chi_\alpha$ in each component $\alpha$ of $S_r = Hom(\RR^\times, U(1))$ (resp{.} $S_c = Hom(\CC^\times, U(1))$). One requires $\chi_\alpha$ to be invariant under positive dilations, so that for the real place $r$ we have $\alpha \in\{+,-\}$ with $\chi_+(u) = 1$ and $\chi_-(u) = sgn(u)$ whereas for the complex place $c$ we have $\alpha = N, N \in \ZZ$ with $\chi_N(z) = z^N(z\overline{z})^{-N/2}$.

The associated Gamma functions have been tabulated by Tate [2a]. He uses the same conventions for the Fourier Transform, but his $\rho(c)$ function is defined in such a way that $\Gamma(\chi,s) = \chi(-1)\rho(\chi(u)|u|^s)$. For the sign character at the real place there is a misprint in the value given in [2a] for $\rho(sgn(u)|u|^s)$ which should have an additional minus sign.
$$\Gamma(\chi_+,s) = \gamma_+(s) = \pi^{{1\over2}-s}\,{\Gamma({s\over2})\over\Gamma({1-s\over2})}$$
$$\Gamma(\chi_-,s) = \gamma_-(s) = i\,\pi^{{1\over2}-s}\,{\Gamma({s+1\over2})\over\Gamma({2-s\over2})}$$
$$\Gamma(\chi_N,s) = \gamma_N(s) = i^{|N|}\,(2\pi)^{1-2s}\,{\Gamma({|N|\over2} + s)\over\Gamma({|N|\over2} + 1-s)}$$

In these equations $\Gamma(s)$ is Euler's Gamma function. Its logarithmic derivative $\lambda(s)$ satisfies
$$\lambda(s) = {\Gamma^\prime(s)\over\Gamma(s)} = -\gamma - {1\over s} - \sum_{j\geq1}({1\over{j+s}}-{1\over j})$$
$$\lambda(s+1) = {1\over s} + \lambda(s)$$
and the spectral functions are
$$\eqalign{
H(\chi_+,s) = h_+(s) &= -\log(\pi) + {1\over2}\lambda({s\over2}) + {1\over2}\lambda({1-s\over2})\cr
K(\chi_+,s) = k_+(s) &= -{i\over4}\left(\lambda^\prime({s\over2}) -  \lambda^\prime({1-s\over2})\right)\cr
H(\chi_-,s) = h_-(s) &= -\log(\pi) + {1\over2}\lambda({1+s\over2}) + {1\over2}\lambda({2-s\over2})\cr
K(\chi_-,s) = k_-(s) &= -{i\over4}\left(\lambda^\prime({1+s\over2}) -  \lambda^\prime({2-s\over2})\right)\cr
H(\chi_N,s) = h_N(s) &= -2\log(2\pi) + \lambda({|N|\over2}+s) + \lambda({|N|\over2} + 1-s)\cr
K(\chi_N,s) = k_N(s) &= -i\left(\lambda^\prime({|N|\over2}+s) - \lambda^\prime({|N|\over2} + 1-s)\right)
}$$
The spectral functions on the critical line for $H$ are
$$\eqalign{
h_+({1\over2}+it) &= -\log(\pi) -\gamma - {1\over{1\over4}+t^2} + \sum_{j\geq1}\left({1\over j}-{4j+1\over(2j+{1\over2})^2 +t^2}\right)\cr
h_-({1\over2}+it) &= -\log(\pi) -\gamma - {3\over{9\over4}+t^2} + \sum_{j\geq1}\left({1\over j}-{4j+3\over(2j+{3\over2})^2 +t^2}\right)\cr
h_N({1\over2}+it) &= -2\log(2\pi) -2\gamma - {|N|+1\over({|N|+1\over2})^2+t^2} + 2\sum_{j\geq1}\left({1\over j}-{j+{|N|+1\over2}\over(j+{|N|+1\over2})^2 + t^2}\right)\cr
}$$
From this we see that they take their minimal values at $t=0$, are even and increase steadily to $\infty$ when $|t|\rightarrow\infty$ (at a logarithmic growth that can also be deduced from Stirling's Formula). These minimal values are 
$$\mu_+ = h_+({1\over2}) = -\log(\pi) + \lambda({1\over4}) = -(\log(8\pi) + \gamma) - {\pi\over2}$$
$$\mu_- = h_-({1\over2}) = -\log(\pi) + \lambda({3\over4}) = -(\log(8\pi) + \gamma) + {\pi\over2}$$
$$\mu_N = h_N({1\over2}) = -2\log(2\pi) + 2\lambda({|N|+1\over2})$$
$$\mu_{2N} = -2(\log(8\pi) + \gamma) + 2\sum_{j=1}^N{1\over j-{1\over2}}\quad (N\geq0)$$
$$\mu_{2N+1} = -2(\log(8\pi) + \gamma) + 4\log(2) + 2\sum_{j=1}^N{1\over j}\quad (N\geq0)$$

{\bf Theorem:} The conductor operator is bounded below.

The spectral functions of the commutator operator on the critical line are 
$$\eqalign{
k_+({1\over2}+it) &= -\sum_{j\geq0}{(4j+1)\,2t\over((2j+{1\over2})^2 +t^2)^2}\cr
k_-({1\over2}+it) &= -\sum_{j\geq0}{(4j+3)\,2t\over((2j+{3\over2})^2 +t^2)^2}\cr
k_N({1\over2}+it) &= -\sum_{j\geq0}{(2j+|N|+1)\,2t\over((j+{|N|+1\over2})^2 + t^2)^2}\cr
}$$
From this or from the Stirling Formula one finds that these functions are $O({1\over|t|})$ when $t\rightarrow\infty$, hence bounded. As $|k_{N+2}(t)| \leq |k_N(t)|$ for $N\geq0$ one gets

{\bf Theorem:} The commutator operator is bounded.

{\bf Spectral analysis in the general case}

In [1a,b,c] a method was used which finds its origin in the Explicit Formula hence uses complex integrals and Mellin Inversion. I will proceed in a more direct manner here. Two certainly very well-known lemmas are needed first.

{\bf Lemma A:} For $\varphi(y)$ a Schwartz-Bruhat function on $K_\nu$ its convolution $B(\varphi)$ with the Fourier Transform of $\log(|x|)$ is a continuous function which is $O({1\over|y|})$ when $|y| \rightarrow \infty$.

{\bf Proof:} As $B(\varphi)$ is the Fourier transform of an $L^1$ function it is continuous. I gave in [1a] a formula for the Fourier transform of $\log(|x|)$. For a finite place it implies that $B(\varphi)(y)$ for $|y|$ large enough is identically ${C\over|y|}$ for some constant $C(\varphi)$. At the real place we write it as
$$-B(\varphi)(y) = \int_{|t| \leq 1}(\varphi(y-t) - \varphi(y))\, {dt \over 2|t|} + \int_{|t| > 1}\varphi(y-t)\, {dt \over 2|t|}+\ (\log(2\pi) + \gamma)\cdot\varphi(y)$$
The first integral is bounded by $\sup_{|t|\geq|y|-1}(|\varphi^\prime(t)|)$. For the second integral one first looks at the contribution of the region of integration $|t| > |y|/2$ which clearly gives $O({1\over|y|})$, and then the remaining part is bounded by $|y|\sup_{|t|\geq|y|/2}(|\varphi(t)|)$ hence with rapid decrease as $|y|\rightarrow\infty$. The complex case is similar and left to the reader (recall though that $|y|$ then means $z\overline{z}$ for $y=z\in\CC$).

{\bf Lemma B:} The identity for $s$ in the critical strip
$$\int \psi(x)\,\chi(x)|x|^{s-1} dx = \Gamma(\chi, s) \int \varphi(y)\,\chi^{-1}(y)|y|^{-s} dy$$
where both $\varphi(y)$ and $\psi(x)$ are supposed to be measurable locally integrable functions, defining tempered distributions, and $\psi(x)$ the Fourier transform of $\varphi(y)$ as a distribution, is valid as soon as both integrals make sense as Lebesgue Integrals (that is are absolutely convergent).

{\bf Proof:} This was indicated in [1c]. One just needs to check that the double integral and change of variables trick of Tate's Thesis works. Let $\alpha(y)$ be a Schwartz-Bruhat function with Fourier Transform $\widetilde{\alpha}(x)$. Starting with
$$\int \psi(x)\,\chi(x)|x|^{s-1} dx \cdot \int \alpha(y)\,\chi^{-1}(y)|y|^{-s} dy$$
we convert it by Fubini's Theorem to a double integral and then apply the change of variables $x = uv, y =v$ to get
$$\int\int\psi(uv)\alpha(v)\chi(u)|u|^{s-1}\,dudv$$
Then one evaluates for $u\neq0$
$$\int\psi(uv)\alpha(v)dv = {1\over|u|}\int\varphi({w\over u})\widetilde{\alpha}(w)dw = \int\varphi(v)\widetilde{\alpha}(uv)\,dv$$
so that the double integral becomes
$$\int\int\varphi(v)\widetilde{\alpha}(uv)\chi(u)|u|^{s-1}\,dudv$$
which is similarly seen (working in reverse) to be
$$\int \widetilde{\alpha}(x)\,\chi(x)|x|^{s-1} dx \cdot \int \varphi(y)\,\chi^{-1}(y)|y|^{-s} dy$$
We have only used the definition of the Fourier Transform of a distribution and the fact that all integrals considered were absolutely convergent, so that the manipulations are allowed by Fubini's Theorem.

Let's start from the fundamental identity for a Schwartz-Bruhat function $\alpha(y)$:
$$\int \widetilde{\alpha}(x)\,\chi(x)|x|^{s-1} dx = \Gamma(\chi, s) \int \alpha(y)\,\chi^{-1}(y)|y|^{-s} dy\leqno{(I)}$$
We are allowed by Lemmas A and B to do two things: to take its derivative, and on the other hand to apply it directly to the Fourier pair $\beta(y) = B(\alpha)(y), \widetilde{\beta}(x) = \log(|x|)\widetilde{\alpha}(x)$. The first operation (with $H(\chi,s)$ defined to be the logarithmic derivative of $\Gamma(\chi,s)$) gives
$$H(\chi,s)\int_{K_\nu}\alpha(y)\chi^{-1}(y)|y|^{-s}\,dy$$
$$= \int_{K_\nu}\alpha(y)\log(|y|)\chi^{-1}(y)|y|^{-s}\,dy + \int_{K_\nu}\widetilde{\alpha}(x)\log(|x|){\chi(x)|x|^{s-1}\over\Gamma(\chi,s)}\,dx$$
and the second gives:
$$\int_{K_\nu}\widetilde{\alpha}(x)\log(|x|){\chi(x)|x|^{s-1}\over\Gamma(\chi,s)}\,dx
= \int_{K_\nu}B(\alpha)(y)\chi^{-1}(y)|y|^{-s}\,dy$$
so that we end up with
$$H(\chi,s)\int_{K_\nu}\alpha(y)\chi^{-1}(y)|y|^{-s}\,dy = \int_{K_\nu}H(\alpha)(y)\chi^{-1}(y)|y|^{-s}\,dy\leqno{(II)}$$
Furthermore by Lemma A, the following integral
$$\int_{K_\nu}H(\alpha_1)(y)\alpha_2(y)\,dy$$
for two test-functions $\alpha_1$ and $\alpha_2$ is absolutely convergent and easily seen to be
$$\int_{K_\nu}\alpha_1(y)H(\alpha_2)(y)\,dy$$
so that the Identity $(II)$ can be interpreted as an identity of distributions
$$H({\chi(y)}^{-1}|y|^{-s}) = H(\chi, s)\ {\chi(y)}^{-1}|y|^{-s}
 = \left({\partial\over\partial s}\log(\Gamma(\chi,s))\right)\ {\chi(y)}^{-1}|y|^{-s}\leqno ({\cal S}_H)$$
This gives a somewhat simpler approach to the spectral analysis of the conductor operator than the one used in [1b].

To proceed further we now restrict $\alpha(y)$ to vanish identically in a neighborhood of the origin so that $A(\alpha)(y) = \log(|y|)\alpha(y)$ is again a Schwartz function. Again on one hand we compute the derivative of $(II)$
$$H^\prime(\chi,s)\int_{K_\nu}\alpha(y)\chi^{-1}(y)|y|^{-s}\,dy$$
$$= H(\chi,s)\int_{K_\nu}\alpha(y)\log(|y|)\chi^{-1}(y)|y|^{-s}\,dy -
\int_{K_\nu}H(\alpha)(y)\log(|y|)\chi^{-1}(y)|y|^{-s}\,dy$$
and on the other hand we apply it directly to $A(\alpha)(y)$:
$$H(\chi,s)\int_{K_\nu}\alpha(y)\log(|y|)\chi^{-1}(y)|y|^{-s}\,dy
= \int_{K_\nu}H(A(\alpha))(y)\chi^{-1}(y)|y|^{-s}\,dy$$
so that we end up with
$$H^\prime(\chi,s)\int_{K_\nu}\alpha(y)\chi^{-1}(y)|y|^{-s}\,dy = \int_{K_\nu}((HA-AH)(\alpha))(y)\chi^{-1}(y)|y|^{-s}\,dy$$
hence with
$$iH^\prime(\chi,s)\int_{K_\nu}\alpha(y)\chi^{-1}(y)|y|^{-s}\,dy = \int_{K_\nu}(K(\alpha))(y)|y|^{-s}\,dy\leqno{(III)}$$
For $\alpha_1$ and $\alpha_2$ vanishing in a neighborhood of the origin one checks easily
$$\int_{K_\nu}K(\alpha_1)(y)\alpha_2(y)\,dy = - \int_{K_\nu}\alpha_1(y)K(\alpha_2)(y)\,dy$$
so that we can reinterpret $(III)$ as an identity of distributions on the {\it punctured} line
$$K({\chi(y)}^{-1}|y|^{-s}) = K(\chi, s)\ {\chi(y)}^{-1}|y|^{-s}
= \left(- i\,{\partial^2\over\partial s^2}\log(\Gamma(\chi,s))\right)\ {\chi(y)}^{-1}|y|^{-s}\leqno ({\cal S}_K)$$
This completes the spectral analysis of the conductor and commutator operators.

{\bf Higher derivatives and higher commutators}

Let's now proceed inductively with higher derivatives and hence deduce the spectral analysis of the higher commutators $K_1 = iK = [A,H], K_{N+1} = [A, K_N]$ as an identity of distributions on the punctured $\nu$-adics
$$K_N({\chi(y)}^{-1}|y|^{-s}) = \left({\partial^{N+1}\over\partial s^{N+1}}\log(\Gamma(\chi,s))\right)\ {\chi(y)}^{-1}|y|^{-s}\leqno ({\cal S}_N)$$
from which the Hilbertian version follows by restriction to the critical line (note that these operators are indeed completely dilation invariant as any commutator built with the identity vanishes).

We apply the operators $K_N$ only to Schwartz functions vanishing identically in a neighborhood of the origin. The explicit formula
$$K_N = \sum_{N\geq j\geq0}(-1)^{N-j}{N\choose j}A^jHA^{N-j}$$
shows that the resulting function is (according to Lemma A) continuous except perhaps at the origin where it is $O(\log(|y|)^N)$ while at infinity it is $O({\log(|y|)^N\over|y|})$ so that Lemma B allows to make sense of $K_N({\chi(y)}^{-1}|y|^{-s})$ as a distribution on the punctured $\nu-$adics and to generalize the method used for $H$ and $K$.

Assuming inductively the validity of $(III_N)$
$$\left({\partial^{N+1}\over\partial s^{N+1}}\log(\Gamma(\chi,s))\right)\int_{K_\nu}\alpha(y)\chi^{-1}(y)|y|^{-s}\,dy = (-1)^N\int_{K_\nu}(K_N(\alpha))(y)\chi^{-1}(y)|y|^{-s}\,dy$$
we first compute its derivative
$$\left({\partial^{N+2}\over\partial s^{N+2}}\log(\Gamma(\chi,s))\right)\int_{K_\nu}\alpha(y)\chi^{-1}(y)|y|^{-s}\,dy =$$
$$\left({\partial^{N+1}\over\partial s^{N+1}}\log(\Gamma(\chi,s))\right)\int_{K_\nu}\alpha(y)\log(|y|)\chi^{-1}(y)|y|^{-s}\,dy$$
$$
- (-1)^N\int_{K_\nu}(K_N(\alpha))(y)\log(|y|)|y|^{-s}\,dy$$
and then apply $(III_N)$ to $A(\alpha)$ so that
$$\left({\partial^{N+2}\over\partial s^{N+2}}\log(\Gamma(\chi,s))\right)\int_{K_\nu}\alpha(y)\chi^{-1}(y)|y|^{-s}\,dy
= (-1)^{N+1}\int_{K_\nu}([A,K_N](\alpha))(y)|y|^{-s}\,dy$$
which is $(III_{N+1})$.

As for $\alpha_1$ and $\alpha_2$ two test-functions with support away from the origin 
$$\int_{K_\nu}K_N(\alpha_1)(y)\alpha_2(y)\,dy = (-1)^N \int_{K_\nu}\alpha_1(y)K_N(\alpha_2)(y)\,dy$$
this gives the identity of distributions $({\cal S}_N)$ on the punctured $\nu-$adics.

As we did before for $K$ from the partial fraction expansion of Euler's Gamma function we see that from the Hilbert Space point of view all $K_N$'s are bounded, and in fact we can say a little more:

{\bf Theorem:} There is a dense subdomain of $L^2(K_\nu, dy)$, stable under the Fourier Transform, $A$, $B$, $H$, $K$, and all higher $K_N$'s. All commutators on this domain of $A$ and $B$ boil down to the already constructed $K_N$'s.

{\bf Proof:} We define it to be the domain of functions whose unitary dual in $L^2(S_\nu, d\mu)$ vanishes except on finitely many components and is of the Schwartz class in each component (which means smooth for a finite place as the components are then circles). Concretely: for a finite place this means, with $q$ the cardinality of the residue field, that we take functions $f(y) = |y|^{-1/2}F({\log(|y|)\over\log(q)})$ where the $F(n), n\in\ZZ$ are the Fourier coefficients of a smooth function on the circle, that is decrease  faster than any polynomial in $|n|^{-1}$ and we also allow twists with multiplicative unitary characters, and then take a finite linear combination. For the real or complex place this means $f(y) = |y|^{-1/2}F(\log(|y|))$ or a twist with a unitary character, and finite linear combinations, where $F(t)$ is a Schwartz function on $\RR$. Clearly this is stable under multiplication with $\log(|y|)$. And this is also stable under additive convolution with the Fourier Transform of $\log(|x|)$ ! Indeed it is enough to show that it is stable under the (additive) Fourier Transform, but in [1b] the Fourier Transform was seen to act in $L^2(S_\nu, d\mu)$ through a composition of Inversion ($\chi \mapsto \overline{\chi}$) and multiplication with the Gamma function (seen on the critical line). For a finite place the Gamma function is smooth on each component so we are done, for an archimedean place all derivatives of $\Gamma(\chi,{1\over2}+it)$ have a growth bounded by a polynomial in $\log(|t|)$. This is seen inductively from $\Gamma^\prime(\chi,s) = H(\chi,s)\Gamma(\chi,s)$ as the $H$ term has logarithmic growth and all its derivatives are bounded (as is easily deduced from the partial fraction expansion of Euler's Gamma function, as seen before), whereas $\Gamma(\chi,s)$ itself is (of course) of absolute value $1$. So multiplying with this leaves the Schwartz space invariant.

It only remains to show that no new commutator can be built. But obviously $[H,K_N] = 0$ as they have the same (generalized) eigenvectors. So any commutator built with $A$'s and $B$'s can be expressed in terms of $A$'s and only one $B$, that is it is one of the $K_N$'s.

{\bf Looking back and making it all (too) easy}

If the goal is just to understand $H$ and $K$ in the Hilbert space, there is a third method which is forced upon us by the proof of the previous theorem.

For $\chi$ a fixed unitary character we parametrize the corresponding component of $S_\nu$ through the association of the character $\chi(u)|u|^{i\tau}$ to the real number $\tau$, taken modulo ${2\pi i \over\log(q)}$ if the place is finite. An element $f(u)$ of $L^2(K_\nu^\times, d^*u)$ has a multiplicative Fourier Transform $f(\chi,\tau) = \int_{K_\nu^\times}f(u)\chi(u)|u|^{i\tau}\,d^*u$. If we restrict to $f$'s in the domain defined above, we can then take the derivative with respect to $\tau$ and end up with:

{\bf Theorem:} In this picture $A$ is ${1\over i}{\partial\over\partial\tau}$

We know from [1b] how the additive Fourier Transform acts:
$$f(\chi,\tau)\mapsto \Gamma(\chi,{1\over2}+i\tau)f(\overline{\chi},-\tau)$$
The inverse Fourier transform acts as
$$f(\chi,\tau)\mapsto \Gamma(\chi,{1\over2}+i\tau)\chi(-1)f(\overline{\chi},-\tau)$$
so $B = {\cal F}A{\cal F}^{-1}$ maps $f(\chi,\tau)$ to
$$\Gamma(\chi,{1\over2}+i\tau)\Gamma^\prime(\overline{\chi},{1\over2}-i\tau)\chi(-1)f(\chi,\tau)
+\Gamma(\chi,{1\over2}+i\tau)\Gamma(\overline{\chi},{1\over2}-i\tau)\chi(-1){1\over i}(-1){\partial\over\partial\tau}f(\chi,\tau)$$
This can be simplified using $\Gamma(\chi,s)\Gamma(\chi^{-1},1-s)=\chi(-1)$ to 
$$H(\chi,{1\over2}+i\tau)f(\chi,\tau) - {1\over i}{\partial\over\partial\tau}f(\chi,\tau)$$
where the function $H(\chi,{1\over2}+i\tau)$ is defined to be the logarithmic derivative of the Gamma function, restricted to the critical line.

So the spectral function of $H = A+B$ is as was previously found $H(\chi,{1\over2}+i\tau)$, and using the Theorem we also get all commutators immediately.

{{\bf REFERENCES}\par
\baselineskip = 12 pt
\parskip = 4 pt
\font\smallRoman = cmr8
\smallRoman
\font\smallBold = cmbx8
\font\smallSlanted = cmsl8
{\smallBold [1a] J.F. Burnol}, {\smallSlanted ``The Explicit Formula and a Propagator''}, math/9809119 (September 1998, revised November 1998).\par
{\smallBold [1b] J.F. Burnol}, {\smallSlanted ``Spectral analysis of the local conductor operator''}, math/9811040 (November 1998).\par
{\smallBold [1c] J.F. Burnol}, {\smallSlanted ``The Explicit Formula in simple terms''}, math/9810169 (October 1998, revised November 1998).\par
{\smallBold [2a] J. Tate}, Thesis, Princeton 1950, reprinted in Algebraic Number Theory, ed. J.W.S. Cassels and A. Fr\"ohlich, Academic Press, (1967).\par
{\smallBold [2b] I. M. Gel'fand, M. I. Graev, I. I. Piateskii-Shapiro},{\smallSlanted ``Representation Theory and automorphic functions''}, Philadelphia, Saunders (1969).\par
{\smallBold [2c] A. Weil},{\smallSlanted ``Fonctions z\^etas et distributions''}, S\'eminaire Bourbaki n${}^{\smallRoman o}$ 312, (1966).\par
{\smallBold [3a] D. Slepian, H. Pollak},{\smallSlanted ``Prolate spheroidal wave functions, Fourier analysis and uncertainty I''}, Bell Syst. Tech. J. 40, (1961).\par
{\smallBold [3b] H.J. Landau, H. Pollak},{\smallSlanted ``Prolate spheroidal wave functions, Fourier analysis and uncertainty II''}, Bell Syst. Tech. J. 40, (1961).\par
{\smallBold [3c] H.J. Landau, H. Pollak},{\smallSlanted ``Prolate spheroidal wave functions, Fourier analysis and uncertainty III''}, Bell Syst. Tech. J. 41, (1962).\par
{\smallBold [3d] A. Connes}, {\smallSlanted ``Trace formula in non-commutative Geometry and the zeros of the Riemann zeta function''}, math/9811068 (November 1998).\par
{\smallBold [3e] H. Dym, H.P. McKean}, {\smallSlanted ``Fourier series and integrals''}, Academic Press (1972).\par
\vfill
\centerline{Jean-Fran\c{c}ois Burnol}
\centerline{62 rue Albert Joly}
\centerline{F-78000 Versailles}
\centerline{France}
\centerline{jf.burnol@dial.oleane.com}
\centerline{December 1998}
}
\eject
\bye